\input amstex
\documentstyle{amsppt}

\magnification\magstep1
\document
\topmatter
\title SOME RATIONAL VERTEX ALGEBRAS
\endtitle

\author Dra\v{z}en Adamovi\'{c} 
\endauthor

\affil
Department of Mathematics, University of Zagreb, Zagreb, Croatia
\endaffil

\address
Department of Mathematics, University of Zagreb, Bijeni\v{c}ka 30, 
41000 Zagreb, Croatia
\endaddress
\email
adamovic\@cromath.math.hr
\endemail



\abstract
Let $L((n-\tfrac 3 2)\Lambda_0)$, $n \in \Bbb N$, be a vertex operator algebra
associated to the irreducible highest weight module $L((n-\tfrac 3 2)\Lambda_0)$
for a symplectic affine Lie algebra.
We find a complete set of irreducible modules for $L((n-\tfrac 3 2)\Lambda_0)$ and
show that every module for $L((n-\tfrac 3 2)\Lambda_0)$ from the category 
$\Cal O$ is  completely reducible.
\endabstract

\endtopmatter
\document
\heading 0. Introduction
\endheading

Let $\goth g$ be a type one affine Lie algebra.
Then the irreducible highest weight \ $\goth g$--module
$L(k\Lambda_0)$ has a natural vertex operator algebra structure 
for every $ k \in \Bbb C$,\ $k \ne -g$.
When  $k$ is a positive integer, then
the vertex operator algebra $L(k\Lambda_0)$ is rational
 and its irreducible modules are integrable highest weight modules of level $k$
(cf. \cite {DL}, \cite {MP}, \cite {FZ}).

In this paper we will consider the case of  a symplectic affine Lie algebra
of the type  $C_{\ell}^{(1)}$ and the corresponding vertex operator algebra
 $L((n-\tfrac 3 2)\Lambda_0)$, $n \in \Bbb N$.
We give the description of  two sets of admissible weights $S_1^n$ and $S_2^n$
and prove that $L(\lambda)$, $\lambda \in S_1^n \cup S_2^n$, are  irreducible
$L((n-\tfrac 3 2)\Lambda_0)$--modules (cf. Section 2 and 3).
Next, we prove that irreducible  $L((n-\tfrac 3 2)\Lambda_0)$--modules
are in one-to-one correspondence with  zeros of the set of  polynomials
$\Cal P_{0,\ell}$ (Section 4).
By using this correspondence we show that the set
$\{L(\lambda)\ \vert \  \lambda \in S_1^n \cup S_2^n \} $ gives a complete list
of irreducible $L((n-\tfrac 3 2)\Lambda_0)$--modules. The classification of
irreducible $L((n-\tfrac 3 2)\Lambda_0)$--modules implies that every 
$L((n-\tfrac 3 2)\Lambda_0)$--module from the category $\Cal O$ is 
completely reducible (cf. Section 6).

It turns out that representations of the vertex operator algebra 
$L((n-\tfrac 3 2)\Lambda_0)$ are in some respects quite "similar" 
to the integrable highest weight representations.  

The author expresses the thanks to M.Primc for suggesting to him the study of
this problem and 
for helpful discussions and valuble comments. 
\heading
1. Symplectic affine algebra
\endheading

The symplectic affine (Kac-Moody) Lie algebra $C^{(1)}_\ell$ can be 
written as

$$
\goth g=sp_{2\ell}(\Bbb C)\otimes\Bbb C[t,t^{-1}] 
+\Bbb C    c +\Bbb C   d $$
with the usual commutation relations (cf. \cite {K}).
For $X \in sp_{2\ell}(\Bbb C) $ and $ n \in \Bbb Z$
we write $X(n)=X \otimes t^n$.
Consider two $\ell$-dimensional vector spaces $A_1=\sum _{i=1}^{\ell} \Bbb C a_i $, $A_2=
\sum _{i=1}^{\ell}\Bbb C    a_i^* $.
Let\  $A=A_1 + A_2$.
The Weyl algebra $W(A)$ is the associative algebra over $\Bbb C$ generated by $A$
 and relations
$$[a_i,a_j]=[a_i^*,a_j^*]=0 ,\quad  [a_i,a_i^*]=\delta_{i,j} ,
\qquad i,j \in \ \{1,2, \dots ,\ell \}. $$
Define the normal ordering on $A$ by 
$$:\!xy\!:\ = \tfrac  1 2 (xy+yx)\ \qquad x,y \in A .$$
Then (cf. \cite {B} and  \cite {FF}) all such elements $:xy:$ span a Lie
algebra isomorphic to $\overset\circ \to { \goth g} =sp_{2\ell}(\Bbb C)$
\ with a Cartan subalgebra\ $\overset \circ \to {\goth h} $\ spanned by
$$h_i=-:a_ia_i^*: \quad   i=1,2,...,\ell.$$
Let $\{\epsilon_i\ | \ 1\leqslant i\leqslant\ell\}\subset
\overset \circ \to {\goth h}^* $ be the dual basis such that
$\epsilon_i(h_j)=\delta_{i,j}$.
The root system of $\overset \circ \to {\goth g}
$ is given by
$$\Delta=\{\pm(\epsilon_i\pm \epsilon_j),\pm 2\epsilon_i\ 
\vert \ 1\leqslant i,j \leqslant
\ell,i<j\}$$
with\ $ \alpha_1=\epsilon_1-\epsilon_2,...,\alpha_{\ell-1}=\epsilon_{\ell-1}-
\epsilon_{\ell},\ \alpha_{\ell}=2\epsilon_{\ell}$
being a set of simple roots. The highest root is $\theta=2\epsilon_1$.
Let $\overset \circ \to {\goth g}=\overset \circ \to {\goth n}_- +
\overset \circ \to {\goth h}+\overset \circ \to {\goth n}_+$
be the corresponding triangular decomposition.
We fix the root vectors :
$$X_{\epsilon_i-\epsilon_j}=\ :\!a_ia_j^*\!:,\quad
X_{\epsilon_i+\epsilon_j}=\ :\!a_ia_j\!:,\quad
X_{-(\epsilon_i+\epsilon_j)}=\ :\!a_i^*a_j^*\!:\ .$$

\heading
{2. Some admissible weights}
\endheading

Let $R$\ (resp $R_+$)$\subset \goth h$ be the set of real\ (resp positive real)
 coroots of\ $ \goth g$. Fix $\lambda \in \goth h^*$. Let 
$R^{\lambda}=\{\alpha\in R\ |\ \langle\lambda,\alpha\rangle\in\Bbb Z\}$, 
$R_+^{\lambda}=R^{\lambda}\cap R_+$,
$\Pi$ the set of simple coroots in $R$  and $\Pi^{\lambda }=
\{\alpha \in R_+^{\lambda}\ \vert \ \alpha$
not equal to a sum of several roots from $R_+^{\lambda}\}$.
Define $\rho$ in the usual way. 

Recall that a weight $\lambda \in \goth h ^* $
is called admissible (cf. \cite {KW 2}) if the following
properties are satisfied : 
\roster
\item $ \langle\lambda + \rho,\alpha\rangle \notin -\Bbb Z_+$ for all 
$\alpha \in R_+$,
\item $\Bbb Q R^{\lambda}=\Bbb Q \Pi $.
\endroster

Let $M(\lambda)$ denote the Verma module with the highest weight $\lambda$, 
 $M^1(\lambda) $ its maximal submodule and $L(\lambda)$ its irreducible quotient.

 First let us recall some results of V.Kac and M.Wakimoto that we shall use:

\proclaim
{Theorem 1}(Kac -- Wakimoto, Cor. 2.1. in \cite {KW 1} )
Let $\lambda $ be an admissible weight. Then
$$L(\lambda)=\frac {M(\lambda)}{\sum_{\alpha \in \Pi^{\lambda}}U (
\goth g) v^{\alpha}}\ ,$$
where $v^{\alpha}\in M(\lambda)$ is a singular vector of weight $r_{\alpha}.
\lambda$, the highest weight vector of $M(r_{\alpha}.\lambda)=\ U(\goth g)
v^{\alpha}\subset M(\lambda)$.
\endproclaim

\proclaim
{Theorem 2}
(Kac--Wakimoto, Theorem 4.1 in \cite {KW 2})
Let $V$ be a $\goth g$-- modul from the category $\Cal O$ such that for any
irreducible subquotient $L(\mu )$ the weight  $\mu $ is admissible. 
Then  $\goth g$--modul $V$ is  completely reducible.
\endproclaim
 Denote by $P_+$ the set of all dominant integral weights. 

\proclaim
{Theorem 3}(Kac--Wakimoto, Cor.4.1 in \cite {KW 2})
Let  $\Lambda \in P_+$ and $\lambda$ be an admissible weight. Then the $\goth g$
--modul $L(\Lambda) \otimes L(\lambda)$ decomposes into a dirtect sum of 
irreducible $\goth g$-- modules $L(\mu)$ with $\mu$ admissible highest weight
and $R^{\mu} = R^{\lambda}$.
\endproclaim

Let $\Pi=\{\alpha_0^{\vee},\dots,\alpha_{\ell}^{\vee}\}$,
$c=\alpha_0^{\vee}+\cdots +\alpha_{\ell}^{\vee}$ and set
$$\Pi_1=\{2c-(h_1+h_2),h_1-h_2,\dots ,h_{\ell-1}-h_{\ell},h_{\ell}\}=
\{2\alpha_0^{\vee}+\alpha_1^{\vee},\alpha_1^{\vee},\dots ,\alpha_{\ell}^{\vee}\},$$
$$\Pi_2=\{c-h_1,h_1-h_2,\dots ,h_{\ell-1}-h_{\ell},h_{\ell-1}+h_{\ell}\}=
\{\alpha_0^{\vee},\dots ,\alpha_{\ell-1}^{\vee},\alpha_{\ell -1}^{\vee}
+2\alpha_{\ell}^{\vee}\}.$$
Let $S_i$ denote the set of all admissible $\lambda$ with 
$\Pi^{\lambda}=\Pi_i, \  i=1,2.$
\proclaim
{Lemma 4}
Let $\lambda \in S_i$, $i=1,2$. Then  
$$\langle \lambda ,c\rangle=n-\tfrac 3 2 \quad \text{for some}\  n \in \Bbb N\ .$$
\endproclaim
\demo {Proof}
 For  $\lambda \in S_1 $ we have  
$$\gathered \langle \lambda +\rho , 2\alpha_0^{\vee}+\alpha_1^{\vee} \rangle =
  \langle \lambda  , 2\alpha_0^{\vee}+\alpha_1^{\vee} \rangle +3 \\
= \langle \lambda  , 2\alpha_0^{\vee}+2\alpha_1^{\vee}+\cdots 
2\alpha_{\ell}^{\vee} \rangle -
\langle \lambda  , \alpha_0^{\vee}+2\alpha_1^{\vee}+\cdots 2\alpha_{\ell}^{\vee} \rangle
+3 >0 .\endgathered  $$
This implies  
$\langle \lambda , c \rangle > -\frac 3 2$ and we see that
 $\langle \lambda , c \rangle \in -\frac 3 2 + \Bbb N$. 
Similary we prove the case i=2. \qed
\enddemo
Let$$ S^n_i=\{\lambda \in S_i \ \vert \ \langle \lambda ,c \rangle = n- \tfrac 3 2 \}
\quad i=1,2,$$ 
$$P_+^1=\{\lambda \in P_+\ \vert \  \langle \lambda ,c \rangle = 1 \}=
\{\Lambda_0,\dots ,\Lambda_{\ell} \}.$$
Then
$S_i=\cup _{n \in \Bbb N}\  S_i^n \ .$
We give a description of  $S_1^n$ and $S_2^n$ for $n \in \Bbb N$ :
\proclaim
{Proposition 5}
 $$
\gathered
 S_1^1=\{-\tfrac 1 2 \Lambda_0,-\tfrac 3 2 \Lambda_0 +\Lambda_1\}\ , \\
S_1^{n+1}=\{S_1^n +P_+^1\}\cup \{-(n+\tfrac 3 2)\Lambda_0 + (2n + 1)  
 \Lambda_1 \},\ \ n \in \Bbb N;
\endgathered \tag1$$
$$
\gathered
S_2^1 =\{-\tfrac 1 2 \Lambda_{\ell},-\tfrac 3 2 \Lambda _{\ell} +\Lambda_{
\ell-1} \}\ , \\
S_2^{n+1}=\{S_2^n + P_+^1\} \cup \{-(n+ \tfrac 3 2)\Lambda_{\ell}+(2n+1)
\Lambda_{\ell-1}\}\ , \quad n \in \Bbb N .
\endgathered \tag2$$

\endproclaim
\demo 
{Proof} We can directly  obtain the description of the  set $S_1^1$.

By the definition of sets $S_i^n$ we have

$$\{S_1^n +P_+^1\} \subset S_1^{n+1}\ \text{and}\  (n+\tfrac 3 2)\Lambda_0 + (2n + 1)  
 \Lambda_1 \in S_1^{n+1}.$$
Let $\lambda \in S_1^{n+1}$, $\lambda \ne -(n+\tfrac 3 2)\Lambda_0 + (2n + 1)  
 \Lambda_1 $. Then 
$\langle \lambda ,\alpha_0^{\vee} \rangle = -(n - m + \tfrac 1 2)$ , for $m \in
\Bbb Z_+$. Since $\langle \lambda + \rho ,
2 \alpha_0^{\vee}+\alpha_1^{\vee} \rangle > 0$ we have
$\langle \lambda ,\alpha_1^{\vee} \rangle \geqslant (2(n-m) -1)$, and this  implies
$$\lambda = -(n-m + \tfrac 1 2)\Lambda_0 + (2(n-m)-1)\Lambda_1 + 
\Lambda^{(1)} + \dots + \Lambda^{(m+1)}$$
where $ \Lambda^{(i)} \in P_+^1,\quad i=1,\dots ,m+1$. We have obtained 
$$\lambda \in S_1^{(n-m)} + P_+^1 + \dots + P_+^1 \subset S_1^n + P_+^1 $$
 and (1) holds.

The proof of (2) is similar. \qed
\enddemo

\heading
{3. Modules for Vertex operator algebra $L((n-\tfrac 3 2)\Lambda_0)$}
\endheading

We know that the generalized Verma module $N(k\Lambda _0)$ 
 with the highest weight
$k\Lambda _0$, $k \in \Bbb C$, is a vertex operator algebra if $k \ne -g$
(here $g$ denotes the dual Coxeter number).
The irreducible quotient $L(k\Lambda _0)$\ of $N(k\Lambda_0)$ is also a vertex
operator algebra (see \cite {FLM}, \cite {FgF},\cite {DL}, 
\cite {FZ} and \cite {MP}). 

As usual we shall denote by $Y(w,z)=\sum_{m \in \Bbb Z}w_m z^{-m-1}$ the 
vertex operator (or the field) of the vector $w$.

Let $V$ be a $\goth g$--module of level $k$, $k \ne -g$ from the category
$\Cal O$ (or a highest weight module) and let 
$$X(z) = Y(X(-1) \bold 1,z)=\sum_{m \in \Bbb Z}X(m) z^{-m-1},\quad  
X \in \overset \circ \to  {\goth g},$$
 be the family of fields acting on $V$
defined by the action of $X(m) \in \goth g$. By Theorem 4.3 in \cite {MP}
or Theorem 2.4.1 in \cite {FZ} there is
a unique extension of these fields that make $V$ into a module over the vertex
operator algebra $N(k\Lambda_0)$. Hence we may identify $\goth g$--modules of level
 $k$ in the category $\Cal O$ with the $N(k\Lambda_0)$--modules in the category
$\Cal O$.

Moreover, if $I$ is an ideal of the vertex operator algebra $N(k\Lambda_0)$, then
a $\goth g$--module from the category $\Cal O$ is a module of the vertex operator
algebra $N(k\Lambda_0)/I$ if and only if $Y(w,z)V=0$ for all $w \in I$ (or 
equivalently, for all generators $w$ of the ideal $I$) (cf. Corrollary 3.2 and
Proposition 4.2 below).

We will find
all irreducible representations of the vertex operator algebras
 $L((n-\tfrac 3 2)\Lambda _0)$, $n \in \Bbb N$,  associated to the symplectic
algebra $C_{\ell}^{(1)}$.

Put $\lambda _n =(n-\tfrac 3 2)\Lambda _0$.\ Then $\lambda _n$ 
is an admissible weight with
$\Pi^{\lambda _n}=\Pi _1 $.

Put $\gamma_0=\delta - (\epsilon_1+\epsilon_2)$.
It is easy to show  that $\gamma_0^{\vee}=2\alpha_0^{\vee}+\alpha_1^{\vee}$.
Then we have:
$$ r_{\gamma_0}.\lambda_n=\lambda_n -2n\gamma_0,\quad
 r_{\alpha_i}.\lambda_n=\lambda_n-\alpha_i,\ \ i=1,2,...,\ell.$$
By $\bold 1$ we denote a highest weight vector in $N(\lambda_n)$.
\proclaim
{Theorem 1}
The maximal submodule of $N(\lambda_n)$ is $N^1 (\lambda_n)=U(\goth g)v_n$,  where
$$v_n=(X_{\epsilon_1+\epsilon_2}(-1)^2 - X_{2\epsilon_1}(-1)X_{2\epsilon_2}(-1)
)^n\bold 1,\ \  \ n \in \Bbb N .$$
\endproclaim
\demo{Proof}
It can be checked by a direct calculation that $v_n$ is a singular vector
of weight \ $ \lambda_n - 2n\gamma_0$. Since 
$$v^{\alpha_i}=X_{-\alpha_i}(0)\bold 1 =0$$
for $i=1,2,...,\ell$, we conclude from Theorem 2.1
that $v_n$ generates the maximal submodule
$N^1(\lambda_n)$.\qed
\enddemo

Clearly we have 
$$ Y(v_n,z) = (X_{\epsilon _1 +\epsilon_2}(z)^2 - X_{2\epsilon_1}(z)
\,X_{2\epsilon_2}(z))^n. $$
Theorem 1 implies the following :
\proclaim
{Corollary 2}
Let $V$ be  $\goth g $--module from the category $\Cal O$ of level $n - \tfrac 3 2 $.
Then
$$( X_{\epsilon _1 +\epsilon_2}(z)^2 - X_{2\epsilon_1}(z)
\,X_{2\epsilon_2}(z))^n=0\quad on \ V $$
if and only if $V$ is $L(\lambda_n)$--module.
\endproclaim

 A.Feingold and I.Frenkel gave
the bosonic construction (see \cite {FF}) of four irreducible 
$\goth g$--modules of level $-\frac 1 2$\
 : $L(\mu_1),L(\mu_2),L(\mu_3),L(\mu_4)$
where 
$$\mu_1=-\tfrac 1 2 \Lambda_0, \ \mu_2=-\tfrac 3 2 \Lambda_0 +\Lambda_1,\ 
\mu_3=-\tfrac 1 2 \Lambda_{\ell},\ \mu_4=-\tfrac 3 2 \Lambda_{\ell}+
\Lambda_{\ell-1}.$$

By using Lemma 7 in \cite{FF} and the explicit construction 
(Theorem A in \cite {FF}) we obtain: 
\proclaim 
{Proposition 3} On $L(\mu_i)$,\quad
$i=1,2,3,4$,
we have
$$ X_{\epsilon _1  +\epsilon_2}(z)^2 - X_{2\epsilon_1}(z)
\,X_{2\epsilon_2}(z)=0\ . $$

\endproclaim

\proclaim
{Theorem 4}
Let $V(n-\tfrac 3 2)$ be an irreducible  $L(\lambda_n)$--module and $V(1)$ an
 irreducible  $L(\Lambda_0)$--module. Then 
$$V(n-\tfrac 3 2)\otimes V(1)$$
is a $L(\lambda_{n+1})$--module.
\endproclaim

\demo 
{Proof}By Theorem 2.2 vector $\bold 1\otimes \bold 1 \in
 L(\lambda_n) \otimes L(\Lambda_0)$ 
generates the submodule isomorphic to $L(\lambda_{n+1})$. It is easy to show that
$L(\lambda_{n+1})$ is a subalgebra of the vertex operator algebra 
$L(\lambda_n) \otimes L(\Lambda_0) $ in the sense of \cite {FZ}. Since
$V(n-\tfrac 3 2)\otimes V(1)$ is a module for $L(\lambda_n) \otimes L(\Lambda_0) $
(cf. Proposition 10.1 in \cite {DL}) it is also a module for $L(\lambda_{n+1})$.\qed  
\enddemo

\proclaim {Remark}
The Theorem 4 can also be proved by using Corrolary 2 
 and the vertex operator formula for integrable highest weight
representations (cf. \cite {LP}, Proposition 5.5).
\endproclaim 
\proclaim
{Lemma 5}
Let $\lambda \in S_1^n \cup S_2^n$. Then\ $L(\lambda)$ is a 
$L(\lambda_n)$--modul.
\endproclaim
\demo
{Proof}
Induction on $n \in \Bbb N$. For $n=1$ we have  $S_1^1 \cup S_2^1 = \{\mu_1,
\mu_2,\mu_3,\mu_4\}$. 
Then $L(\mu_i)$, $i=1,2,3,4$ are $L(-\tfrac 1 2 \Lambda_0)$ modules by
Proposition 3.

First notice that  $L(\Lambda)$
 for $\Lambda \in P_+^1$ is a $L(\Lambda_0)$--module
(cf. \cite {FZ}, \cite {DL}, \cite {MP}). 
Assume that $L(\lambda ')$ is a $L(\lambda_n)$--modul for all $\lambda'
\in S_1^n \cup S_2^n$.
Let $\lambda \in S_1^{n+1} \cup S_2^{n+1}$.
If $\lambda = \lambda_0 +\Lambda $, $\lambda_0 \in S_1^n \cup S_2^n$,
 $\Lambda \in P_+^1$, then $L(\lambda_0)\ \otimes L(\Lambda)$ 
is a $L(\lambda_{n+1})$--module
by Theorem 4. Since $v_{\lambda_0} \otimes v_{\Lambda}$ is a singular vector
of weight $\lambda_0 + \Lambda $, by Theorem 2.2 it generates the submodule isomorphic to
 $L(\lambda)$ and $L(\lambda)$ is a $L(\lambda_{n+1})$--module.

Let $\lambda = -(n + \tfrac 3 2)\Lambda_0 + (2n+1)\Lambda_1$.
Put $\mu =-(n + \tfrac 1 2)\Lambda_0 + (2n-1)\Lambda_1$.
Then $L(\mu) \otimes L(\Lambda_0)$ is a  $L(\lambda_{n+1})$--module.
Since 
$$v=2X_{2\epsilon_1}(-1)v_{\mu} \otimes v_{\Lambda_0}+
(2n+1)v_{\mu} \otimes X_{2\epsilon_1}(-1)v_{\Lambda_0}$$
is a singular vector in $L(\mu) \otimes L(\Lambda_0)$ of weight 
$\lambda - \delta $, it generates the submodule isomorphic to $L(\lambda)$.

Let $\lambda=-(n+ \tfrac 3 2)\Lambda_{\ell} +(2n+1)\Lambda_{\ell-1}$.
Put $\nu =-(n+ \tfrac 1 2)\Lambda _{\ell} + (2n-1)\Lambda_{\ell-1}$. Then
$L(\nu) \otimes L(\Lambda_{\ell})$ is a $L(\lambda_{n+1})$--modul.
Since
$$v_1=2X_{-2\epsilon_{\ell}}(0)v_{\nu} \otimes v_{\Lambda_{\ell}}
+(2n+1)v_{\nu} \otimes X_{-2\epsilon_{\ell}}(0)v_{\Lambda_{\ell}}$$
is a singular vector in $L(\nu) \otimes L(\Lambda_{\ell})$, it
generates the submodule isomorphic to $L(\lambda)$.\qed

\enddemo
\proclaim
{Remark} 
It follows from Lemma 5  that $L(\lambda)$, $\lambda \in S_1^n 
\cup S_2^n$, is an irreducible $L(\lambda _n)$--module.
 In what follows we prove that these are all irreducible
 $L(\lambda_n)$--modules (cf. Lemma 6.1 and Theorem 6.2).
\endproclaim

\heading
4. Classification of irreducible representations
\endheading

Fix $n \in \Bbb N$. For  $w \in U(\overset \circ\to {\goth g})v_n$ and
$j \in \Bbb Z$ put $w(j)=w_{j+2n-1}$. Then $w(j)$  has operator degree $j$
(i.e. $[d,w(j)]=jw(j)$).

Set
$$ \overline W = \coprod _{j \in \Bbb Z} W(j),\quad W(j)=
\Bbb C\! - \! span \ \{  w(j)   \vert \  
w \in U(\overset \circ\to {\goth g}) v_n\ \}.$$

By using the commutator formula for vertex operators we get (cf. \cite {MP})
the following:
\proclaim
{Proposition 1}
$\overline W$\   is a loop module under the adjoint action of $\goth g$.
In particular,
$$[X(i),w(j)]=(X.w)(i+j)$$
for $X \in \overset \circ\to {\goth g}$,
$w \in  U(\overset \circ\to {\goth g}) v_n$, $i,j \in \Bbb Z$.
\endproclaim

Then $W(0)$ is a finite dimensional $\overset \circ \to {\goth g}$--module
with the highest weight \ $2n(\epsilon _1 + \epsilon _2) = 2n \ \omega _2$. 
By $W(0)_0$ denote the zero-weight subspace of $W(0)$.

\proclaim
{Proposition 2}
Let $V$ be an irreducible highest weight module of level $n-\tfrac 3 2 $ with the highest weight
vector $v_{\lambda}$.
 The following statements are equivalent :
\roster
\item $V$ is a $L(\lambda _n)$--module;
\item $\overline W V = 0\ $;
\item $W(0)_0 v_{\lambda } = 0$.
\endroster
\endproclaim
\demo {Proof}

The equivalence of (1) and (2) was already discussed in the introduction of
Section 3.

Clearly (2) implies (3).

 For the converse first notice that by assumption
$V=M(\lambda)/M^1(\lambda)$. Hence to see (2) it is enough to see 
$ \overline W M(\lambda) \subset M^1(\lambda)$, i.e. 
$\overline W M(\lambda) \ne M(\lambda)$ (since 
$\overline W M(\lambda)$ is a submodule, and  $M^1(\lambda)$ is the maximal
submodule).

Since $\overline W M(\lambda)=\overline W U(\goth n_-)v_{\lambda}=
U(\goth n_-)\overline W v_{\lambda}$, we have

$$\overline W M(\lambda)\ne M(\lambda)\quad \text{iff} \quad
 v_{\lambda} \in \overline W M(\lambda) \quad \text{iff} \quad
 W(0)_0v_{\lambda}=0. \qed $$

\enddemo

Let $ u \in W(0)_0 $. Clearly there exists the uniqe polynomial
$p_u \in S(\overset \circ \to {\goth h} )$ such that 
$$uv_{\lambda}\  =\  p_u(\lambda)v_{\lambda}.$$ Set
 $ \Cal P_{0,\ell} = \{\ p_u \ \vert \ u \in W(0)_0 \}$. We have
\proclaim
{Corollary 3}
There is one-to-one correspondence between :
\roster
\item   irreducible $L(\lambda_n)$--modules from the category $\Cal O$;

\item $\lambda \in \goth h^* $ such that
$p(\lambda)=0$ \ \ for all $p \in \Cal P _{0,\ell}$.
\endroster
\endproclaim

\heading
{5.Zeros of some polynomials}
\endheading
Denote by $_L$ the adjoint action of $\overset \circ \to {\goth g}$ on
$U(\overset \circ \to {\goth g})$ : $ X_Lf=[X,f]$ for $X \in 
\overset \circ \to {\goth g}$ and 
$f \in U(\overset \circ \to {\goth g})$.
The following lemma is obtained by direct calculations:
\proclaim {Lemma 1}
\roster
\item $ (X_{2\epsilon_1}^k) _L (X_{-2\epsilon_1}^n) \in 
X_{-2\epsilon_1}^{n-k}(-1)^k4^kn\cdots (n-k+1) 
\cdot (h_1-n+k)\cdots (h_1-n+1) 
+ U(\overset \circ \to {\goth g})\overset \circ \to {\goth n}_+ $,
\item $ (X_{2\epsilon_1}^n) _L (X_{-2\epsilon_1}^n) \in
(-1)^n4^n n!\cdot h_1\cdots (h_1 -n+1) 
+U(\overset \circ \to {\goth g})\overset \circ \to {\goth n}_+$
\item $ (X_{2\epsilon_2}^n) _L (X_{-2\epsilon_2}^n) \in
(-1)^n4^n n!\cdot h_2\cdots (h_2-n+1) 
+U(\overset \circ \to {\goth g})\overset \circ \to {\goth n}_+$,
\item $ (X_{\epsilon_1 +\epsilon_2}^m) _L (X_{-\epsilon_1 -\epsilon_2}^m) \in
(-1)^m m!\cdot (h_1 +h_2 )\cdots (h_1+h_2-m+1) 
+U(\overset \circ \to {\goth g})\overset \circ \to {\goth n}_+$,
\item $ (X_{\epsilon_1 +\epsilon_2}^{m'}) _L (X_{-\epsilon_1 -\epsilon_2}^m) \in
U(\overset \circ \to {\goth g})X_{\epsilon_1 +\epsilon_2}$
for $m'>m$, 
\item $(X_{\epsilon_1 +\epsilon_2}^r)_L (X_{-2\epsilon_2}^k) \in
U(\overset \circ \to {\goth g})\overset \circ \to {\goth n}_+$ for $r>0$,
\item $(X_{\epsilon_1 +\epsilon_2}^{2k})_L (X_{-2\epsilon_1}^k)
=(2k)!X_{2\epsilon_2}^k$,
\item $(X_{\epsilon_1 +\epsilon_2}^{2k+i})_L (X_{-2\epsilon_1}^k)
=0$ for $i>0$,
\item $p(h)X_{\alpha}^k=X_{\alpha}^k p(h+k\alpha(h))$ for any
 polynomial $p$.
\endroster
\endproclaim
\proclaim {Lemma 2}
Let $$\gathered f=X_{\beta_1}\cdots X_{\beta_k},\quad X_{\beta_i}\in\overset \circ \to {\goth n}_+, \quad
[X_{\beta_i},X_{\beta_j}]=0 ,\quad \text{for all}\quad i,j; \\
g=X_{\gamma_1}\cdots X_{\gamma_m},\quad X_{\gamma_i}\in\overset \circ \to {\goth n}_-, \quad
[X_{\gamma_i},X_{\gamma_j}]=0 ,\quad \text{for all}\quad i,j;\\
 \sum_{i=1}^k\beta_i + \sum_{i=1}^m\gamma_i = 0. \endgathered $$
Then
\roster
\item $f_Lg \in X_{\beta_1}\cdots X_{\beta_k}X_{\gamma_1}\cdots X_{\gamma_m}
+U(\overset \circ \to {\goth g})\overset \circ \to {\goth n}_+$,
\item $g_Lf \in (-1)^mX_{\beta_1}\cdots X_{\beta_k}
X_{\gamma_1}\cdots X_{\gamma_m} +U(\overset \circ \to {\goth g})
\overset \circ \to {\goth n}_+$.
\endroster
\endproclaim

We shall also use the following consequence of the binomial formula :
\proclaim {Lemma 3} For a polynomial $q$ of degree  $\deg (q) <n$ we have
$$\sum_{k=0}^n(-1)^k \binom n k q(k) = 0.$$
\endproclaim 

In this Section we consider the case $C_2$ and calculate 
some polynomials from $\Cal P _{0,2}$.

\proclaim
{Lemma 4}
Let:
\roster
\item
$p_1(h)=(h_1-h_2)(h_1-h_2-1) \dots (h_1-h_2-2n+1)$;
\item
$p_2(h)=(h_1-n+\frac 3 2)(h_1-n+\frac 5 2)\dots (h_1+\frac 1 2)
h_2(h_2-1) \dots (h_2-n+1)$;
\item
$p_3(h)=\sum_{k=0}^n \frac {n!4^n}{k!4^k} (h_1+h_2-2n+1)  \dots (h_1+h_2 -2n+2k)
h_2(h_2-1)\dots (h_2-n+k+1)$.
\endroster
Then $p_1,p_2,p_3 \in \Cal P _{0,2}$.
\endproclaim
\demo
{Proof}

 We identify  $\overset \circ \to {\goth g} \otimes t^0$
with $\overset \circ \to {\goth g}$ and write $X$ instead of $X(0)$.
Clearly for $a_1,a_2\dots ,a_r \in \overset \circ \to {\goth g}$ we have 
$$\gathered [\text{Coeff}_{z^{-2n}}\quad  (a_1\cdot a_2 \cdots a_r)_L 
(X_{\epsilon_1 + \epsilon_2}(z) ^2 - X_{2\epsilon_1}(z)X_{2\epsilon_2}(z))^n]
v_{\lambda}\\
=[(a_1\cdot a_2 \cdots a_r)_L
(X_{\epsilon_1 + \epsilon_2}(0)^2 - X_{2\epsilon_1}(0)X_{2\epsilon_2}(0))^n]
v_{\lambda }. \endgathered$$ 
Hence
$$W(0)v_{\lambda}=Wv_{\lambda},\qquad W(0)_0 v_{\lambda}=W_0v_{\lambda},$$
where 
$$W=U(\overset \circ \to {\goth g})_L
(X_{\epsilon_1 + \epsilon_2}^2 - X_{2\epsilon_1}X_{2\epsilon_2})^n
\subset U(\overset \circ \to {\goth g})$$
and where $W_0$ denotes the zero weight subspace of $W$.

(1) First notice that
$$(X_{\epsilon_1 - \epsilon_2} ^{2n} X_{-2\epsilon_1}^{2n})_L 
 (X_{\epsilon_1 + \epsilon_2} ^2 - X_{2\epsilon_1}X_{2\epsilon_2})^n 
=8^n (X_{\epsilon_1 - \epsilon_2} ^{2n})_L
(X_{-\epsilon_1 + \epsilon_2} ^2 - X_{-2\epsilon_1}X_{2\epsilon_2})^n \in W_0 $$
We have 
$$\gathered (X_{\epsilon_1  - \epsilon_2}^{2n}) _L 
(X_{-\epsilon_1 + \epsilon_2} ^2 - X_{-2\epsilon_1}X_{2\epsilon_2})^n \\ = 
 (X_{\epsilon_1 - \epsilon_2}^{2n})_L  X_{-\epsilon_1 + \epsilon_2}^{2n} +
\sum_{k=0}^{n-1} \binom  nk (-1)^k (X_{\epsilon_1 - \epsilon_2}^{2n})_L 
X_{-\epsilon_1 + \epsilon_2}^{2k} X_{-2\epsilon_1}^{n-k} X_{2\epsilon_2}^{n-k}\\
\in Cp_1(h) + U(\overset \circ \to {\goth g})\overset \circ \to {\goth n} _+
 \endgathered $$
for some  constant $C\ne 0$.

(2) First notice that
$$(X_{-2\epsilon_1}^n X_{-2\epsilon_2}^n) _L
 (X_{\epsilon_1 + \epsilon_2} ^2 - X_{2\epsilon_1}X_{2\epsilon_2})^n 
 \in W_0.$$
By  Lemma 2 we may calculate  the corresponding polynomial from
$$\gathered (X_{\epsilon_1 + \epsilon_2} ^2 - X_{2\epsilon_1}X_{2\epsilon_2})^n _L 
(X_{-2\epsilon_1}^n X_{-2\epsilon_2}^n) \\
=\sum_{k=0}^n (-1)^k\binom n k (X_{2\epsilon_2}^kX_{\epsilon_1 + \epsilon_2}^{2n-2k}
X_{2\epsilon_1}^k)_L(X_{-2\epsilon_1}^n X_{-2\epsilon_2}^n). \endgathered$$
By using Lemma 1 we have:
$$\gathered (X_{2\epsilon_1}^k)_LX_{-2\epsilon_1}^n X_{-2\epsilon_2}^n
\in X_{-2\epsilon_1}^{n-k} X_{-2\epsilon_2}^n \\
\cdot (-1)^k4^kn(n-1)\cdots (n-k+1)(h_1-n+k)\cdots (h_1-n+1) +
U(\overset \circ \to {\goth g})\overset \circ \to {\goth n} _+ \endgathered$$
and
$$\gathered
(X_{\epsilon_1 + \epsilon_2}^{2n-2k})_L(X_{-2\epsilon_1}^{n-k} X_{-2\epsilon_2}^n)
=[(X_{\epsilon_1 + \epsilon_2}^{2n-2k})_L(X_{-2\epsilon_1}^{n-k}]
 X_{-2\epsilon_2}^n + \\
\sum_{i=1}^{2n-2k}\binom {2n-2k}{i}
[(X_{\epsilon_1 + \epsilon_2}^{2n-2k-i})_L(X_{-2\epsilon_1}^{n-k}]
[(X_{\epsilon_1 + \epsilon_2}^i)_LX_{-2\epsilon_2}^n]\\
\in  [(X_{\epsilon_1 + \epsilon_2}^{2n-2k})_LX_{-2\epsilon_1}^{n-k}]
 X_{-2\epsilon_2}^n +
U(\overset \circ \to {\goth g})\overset \circ \to {\goth n} _+
\quad (\text{ by Lemma 1.6}) \\
=(2n-2k)!X_{2\epsilon_2}^{n-k}X_{-2\epsilon_2}^n
+U(\overset \circ \to {\goth g})\overset \circ \to {\goth n} _+
\quad (\text{by Lemma 1.7}). \endgathered $$
We have obtained
$$\gathered
(X_{2\epsilon_2}^kX_{\epsilon_1 + \epsilon_2}^{2n-2k}X_{2\epsilon_1}^k)_L
X_{-2\epsilon_1}^n X_{-2\epsilon_2}^n \\
\in X_{2\epsilon_2}^{n}X_{-2\epsilon_2}^n (-1)^k4^kn(n-1)\cdots (n-k+1)\\
\cdot (2n-2k)!(h_1-n+k)\cdots (h_1-n+1) +
U(\overset \circ \to {\goth g})\overset \circ \to {\goth n} _+ \\
= (-1)^n n! 4^n h_2\cdots (h_2 -n+1) \\ \cdot
\sum_{k=0}^n \binom n k (2n-2k)!4^k n\cdots(n-k+1) (h_1-n+k)\cdots(h_1-n+1) +
U(\overset \circ \to {\goth g})\overset \circ \to {\goth n} _+. \endgathered$$

For $ h_1=-\tfrac 1 2 + j$, $j=0,1,\dots ,n$ we can show 
$$\sum_{k=0}^n \binom n k (2n-2k)!4^k n\cdots(n-k+1) (h_1-n+1)\cdots(h_1-n+k) $$
$$ \split
&=(2n)!!\sum_{k=0}^n \binom n k (2n-2k-1)!! (-2n+2+2j-1)\cdots (-2n+2k+2j-1) \\
&=(2n)!!(2n-2j-1)!!\sum_{k=0}^n (-1)^k \binom n k (2n-2k-1)\cdots (2n+1-2j-2k) \\
&=0 \qquad (\text{by using Lemma 3}) . \endsplit $$
This implies 
$$\gathered \sum_{k=0}^n \binom n k (2n-2k)!4^k n\cdots(n-k+1) 
(h_1-n+1)\cdots(h_1-n+k) \\ =
4^n n!(h_1-n+\tfrac 3 2)(h_1-n+\tfrac 5 2)\dots (h_1+\tfrac 1 2) \endgathered $$
and we have 
$$(X_{-2\epsilon_1}^n X_{-2\epsilon_2}^n) _L
 (X_{\epsilon_1 + \epsilon_2} ^2 - X_{2\epsilon_1}X_{2\epsilon_2})^n  
\in C p_2(h) +
U(\overset \circ \to {\goth g})\overset \circ \to {\goth n} _+ $$
for some  constant $C\ne 0$.

(3) First notice that
$$(X_{\epsilon_1 + \epsilon_2}^{2n})_L
(X_{-\epsilon_1 - \epsilon_2} ^2 - X_{-2\epsilon_1}X_{-2\epsilon_2})^n
 \in W_0.$$
By using Lemma 1 we can show
$$\gathered 
(X_{\epsilon_1 + \epsilon_2}^{2n})_L(X_{-2\epsilon_1}^k
X_{-\epsilon_1 - \epsilon_2} ^{2n-2k}X_{-2\epsilon_2}^k) \\
\in (2n)!4^kk!(-1)^kh_2\cdot(h_2-k+1)(h_1+h_2-2n+1)\cdots (h_1+h_2-2k) +
U(\overset \circ \to {\goth g})\overset \circ \to {\goth n} _+ .\endgathered$$
By using this and Lemma 1 we see that 
$$(X_{\epsilon_1 + \epsilon_2}^{2n})_L
(X_{-\epsilon_1 - \epsilon_2} ^2 - X_{-2\epsilon_1}X_{-2\epsilon_2})^n
\in (2n)!p_3(h)+
U(\overset \circ \to {\goth g})\overset \circ \to {\goth n} _+ . \qed$$
\enddemo

The following lemma describes the set
$$T^n=\{h \in \Bbb C ^2\  \ |\ p_1(h)=p_2(h)=p_3(h)=0\}.$$
\proclaim
{Lemma 5}
$T^n=T_1^n \cup T_2^n$, where
$$\gathered T_1^n=\{(s+2r,,s) \  |\ s=0,1, \dots ,n-r-1\ ,\ r=0,1,\dots ,n-1 \}\\
  \cup \quad \{(s+2r+1,s) \ |\ s =0,1, \dots ,n-r-1 \ , \ r=0,1, \dots ,n-1 \},\\
 T_2^n=\{(s+2r,s) \ |\ s=-r-\tfrac 1 2 ,\dots,n-2r- \tfrac 3 2 ,\ 
r=0,1,\dots ,n-1 \}\\
 \cup \quad \{(s+2r+1,s)\  |\ s=-r- \tfrac 3 2 ,\dots ,n-2r-\tfrac 5 2 ,
\ r=0,1,\dots ,n-1\}.\endgathered$$
\endproclaim
\demo {Proof}Fix $n \in \Bbb N$ and
let $T_{1,2}=\{h \in \Bbb C ^2\  \ |\ p_1(h)=p_2(h)=0\}$. Then
 $T_{1,2}=T^1_{1,2} \cup T^2_{1,2}$ where 
$$T^1_{1,2}=\{(k,k') \in \Bbb Z^2 \ \vert \ k'=0,\dots,n-1,
k=k',\dots,k' +2n-1\},$$
$$T^2_{1,2}=\{(k,k') \in \Bbb (Z+\tfrac 1 2)^2 \ \vert \ k=-\tfrac 1 2 +i,
\ i=0,\dots,n-1;\ k'=k-j, \ j=0,\dots,2n-1 \}.$$
Clearly $h \in T^n$ if and only if $p_3(h)=0$.

Let $(h_1,h_2) \in T_{1,2} $ and $h_1-h_2=2r,\quad r=0,\dots,n-1$.
 Put $h_2 = s$. Then for 
$\tilde {p_3}(s)=p_3(s+2r,s)$ we have
$$\gathered \tilde {p_3}(s)=\sum_{k=0}^n \frac {n!4^n}{k!4^k}(2s-2n+2r+1)\cdots 
(2s-2n+2r+2k) \cdot  \cdots\ (s-n+k+1) \\
=4^n(n-r)!r!\binom {s}{n-r} \cdot 
\sum_{k=0}^n \binom {s-n+r+k}{r}
\binom {s-n+r+k-\tfrac 1 2}{k} .\endgathered$$
Let $(2r+s,s) \in T_{1,2}$. Clearly $(2r+s,s) \in T^n$ if and only if 
$\tilde {p_3}(s)=0$.
It is easy to see that
$$ \tilde {p_3}(s)=0\ \text{for}\  s=0,\dots,n-r-1 \quad\text {and}\quad
\tilde {p_3}(s)\ne 0\  \text{for}\  s=n-r,\dots,n-1.$$

Let $s=-r-\tfrac 1 2 +i,\quad i=0,\dots ,n-r-1$. Then we have
$$\gathered \frac {1}{4^n(n-r)!r!}\tilde {p_3}(s)=\binom {s}{n-r}\sum_{k=0}^n 
\binom {-n-\tfrac 1 2 +i+k}{r}\binom {-n-1+i+k}{k} \\
=\binom {s}{n-r}\sum_{k=0}^n (-1)^k\binom{n-i}{k} \binom {-n-\tfrac 1 2 +i+k}
{r} = 0 \endgathered $$
by using  Lemma 3.

Let $s=-r-\tfrac 1 2 -i,\quad i=1,\dots,r$. Then we have (by using Lemma 3)
$$\frac {1}{\binom {s}{n-r}4^n(n-r)!r!}\tilde {p_3}(s)
=\sum_{k=0}^n (-1)^k\binom{n+i}{k}
 \binom {-n-\tfrac 1 2 -i+k}{r}$$
$$\split
&=\sum_{k=0}^{n+i} (-1)^k\binom{n+i}{k}
 \binom {-n-\tfrac 1 2 -i+k}{r}- 
\sum_{k=n+1}^{n+i} (-1)^k\binom{n+i}{k}\binom {-n-\tfrac 1 2 -i+k}{r}\\
&=-\sum_{k=n+1}^{n+i} (-1)^k\binom{n+i}{k}\binom {-n-\tfrac 1 2 -i+k}{r} \\
&=(-1)^{n+i+r+1}\sum_{k=0}^{i-1} (-1)^k\binom{n+i}{k}
\binom {k+r-\tfrac 1 2 }{r}.
\endsplit $$
Since
$$\frac {\binom{n+i}{k+1} \binom {k+r+\tfrac 1 2}{r}}
 {\binom{n+i}{k} \binom {k+r-\tfrac 1 2}{r}}
=\frac{(n+i-k)(k+r+\tfrac 1 2)}{(k+1)(k+\tfrac 1 2)}>1$$
we can easily show that $\tilde {p_3}(s)\ne 0$.

Similarly we treat the case $$(h_1,h_2) \in T_{1,2}, \quad h_1-h_2=2r+1, 
\quad r=0,\dots,n-1$$
and  obtain the result.\qed
\enddemo

It follows from  Lemma 5 :
\proclaim
{Lemma 6}
\roster
\item
$T_1^{n+1}=T_1^n \cup \{T_1^n +(1,0)\} \cup \{T_1^n +(1,1)\}\cup
\{(2n+1,0)\}$;
\item
$T_2^{n+1}=T_2^n \cup \{T_2^n +(1,0)\} \cup \{T_2^n +(1,1)\}\cup
\{(n-\frac 1 2 ,-n - \frac 3 2 )\}$.
\endroster
\endproclaim

\heading
{6.The main result}
\endheading

\proclaim
{Lemma 1}
Let $L(\lambda )$ be a $L(\lambda_n)$--module. Then $\lambda \in S_1^n
\cup S_2^n. $
\endproclaim

\demo {Proof} 
Let $L(\lambda )$ be a $L(\lambda_n)$--module. Since 
$$(X_{\epsilon_j + \epsilon_{j+1}}(0)^2 - X_{2\epsilon_j}(0)X_{2\epsilon_{j+1}}(0))^n
\quad \in W, \quad  j=1,\dots ,\ell $$
we can use results for the case $C_2$ and   obtain that 
$$(\lambda (h_j),\lambda (h_{j+1}) )\in T^n ,\quad 1\leqslant j \leqslant \ell.$$
Let $\tilde S _i ^n =\{\lambda \in \goth h^* \vert \ \langle \lambda,c\rangle
=n-\tfrac 3 2,\ \ 
(\lambda (h_j),\lambda (h_{j+1}) )\in T_i^n ,\quad 1\leqslant j \leqslant \ell \}$,
 i=1,2. We will prove by induction that $\tilde S _i ^n= S _i ^n $ for all $n \in \Bbb N $,
$i=1,2$.

For $n=1$ we have $T_1^1=\{(0,0),(1,0)\}$ and  $T_2^1=\{(- \frac 1 2,- \frac 1 2),
(- \frac 1 2,- \frac 3 2)\}$. Then for $\lambda \in 
\tilde S _1 ^1 \cup \tilde S _2 ^1 $ we get
$$(\lambda(h_1),\dots,\lambda(h_{\ell})) \in 
\{(0,\dots,0),(1,0,\dots,0),(-\tfrac 1 2,\dots,-\tfrac 1 2),(-\tfrac 1 2,\dots,
-\tfrac 1 2,-\tfrac 3 2) \}$$ which implies
$$\gathered \tilde S _1 ^1 =\{-\tfrac 1 2 \Lambda_0,-\tfrac 3 2 \Lambda_0 +
\Lambda_1 \}= S_1^1 , \\
\tilde S _2 ^1 =\{-\tfrac 1 2 \Lambda_{\ell},-\tfrac 3 2 \Lambda_{\ell} +
\Lambda_{\ell-1} \}= S_2^1. \endgathered $$

Assume that $\tilde S _i^n=S_i^n$. 
Let $\lambda \in \tilde S _1^{n+1}$. 

If  $(\lambda(h_1),\lambda(h_2))=(2n+1,0)$ then 
$\lambda = -(n+\tfrac 3 2)\Lambda_0 +(2n+1)\Lambda_1.$

If $(\lambda(h_1),\lambda(h_2))\ne(2n+1,0)$ then
$(\lambda (h_j),\lambda (h_{j+1}) )\in
T_1^n \cup \{T_1^n +(1,0)\} \cup \{T_1^n +(1,1)\}$ $j=1,\dots,\ell$.

 We define $\Lambda \in \goth h ^*$ by
$$\split
\langle \Lambda ,h_j \rangle& =\left\{ \aligned & 0 \quad \text{if}\quad 
(\lambda (h_j),\lambda (h_{j+1})
\in T_1^n \\
& 1 \quad \text {otherwise} \endaligned \right. \quad
\text{for}\  j=1,\dots,\ell-1, \\
\langle \Lambda ,h_{\ell}\rangle & =\left\{ \aligned & 0 \quad \text{if}\quad 
(\lambda (h_{\ell-1}),\lambda (h_{\ell})
\in \{T_1^n +(1,0)\}\cup T_1^n \\
& 1 \quad \text {otherwise} \endaligned \right. ,\\
\langle \Lambda ,c \rangle &=1 .
\endsplit  $$ 
Let $\lambda' = \lambda - \Lambda$. It is easy to show that  $\Lambda \in P_+^1$
and $\lambda' \in\tilde S _1^n $. 
Since $\{\tilde S _1^{n} +P_+^1\} \subset \tilde S _1^{n+1}$  we have obtained 
$$\gathered \tilde S _1^{n+1} = \{\tilde S _1^{n} +P_+^1\} \cup \{-(n+\tfrac 3 2)
\Lambda_0 + (2n+1)\Lambda_1\} \\
=\{ S _1^{n} +P_+^1\} \cup \{-(n+\tfrac 3 2)\Lambda_0 + (2n+1)\Lambda_1\} =
 S_1^{n+1} \endgathered $$ 
(by using Proposition 2.5) .

Similary we prove
$$ \tilde S _2^{n+1} = \{\tilde S _2^{n} +P_+^1\} \cup \{-(n+\tfrac 3 2)
\Lambda_{\ell} + (2n+1)\Lambda_{\ell-1}\}= S_2^{n+1}$$
and we conclude by induction that 
 $\tilde S _i^n=S_i^n$, $i=1,2$.
\qed
\enddemo
\proclaim
{Theorem 2}
\roster
\item
The set $\{L(\lambda)\ \vert \  \lambda \in S_1^n \cup S_2^n \}$
provides a complete list of irreducible $L(\lambda_n)$--modules.
\item
Let $V$\ be a $L(\lambda _n)$--module 
from the category $\Cal O $. Then $V$ decomposes into a direct
sum of irreducible $L(\lambda_n) $--modules.
\endroster
\endproclaim
\demo {Proof}(1) By using Lemma 3.5 and Lemma 1 we have that 
$L(\lambda_n)$--modules are exactly $L(\lambda)$ for 
$\lambda \in S_1^n \cup S_2^n$.

(2) Let $L(\mu)$ be an irreducible subquotient of\  $V$.  Then $L(\mu)$ 
is a $L(\lambda _n)$--module and by Lemma 1 we have that
 $\mu \in S^n_1 \cup S^n_2 $.
By using Theorem 2.2 we obtain that $V$ is  completely reducible.\qed
\enddemo

\proclaim
{Remark}
In \cite {Z} and \cite {FZ} are defined representations of vertex 
operator algebras which need not be in category $\Cal O$.
 Vertex operator algebra is by definition
rational if it has only finitely many irreducible modules and if every finitely
generated module is a direct sum of irreducible ones. By the abuse of language
(or by changing the definition) one could say that Theorem 2 states that 
the vertex operator algebra $L((n-\tfrac 3 2)\Lambda_0)$, $n \in  \Bbb N$,  is
rational.
\endproclaim
By using Theorem 2 we obtain:

\proclaim
{Corollary 3 }
Let $V$ be a highest weight $\goth g$--module of level $n-\tfrac 3 2 $.
 The following statements are equivalent :
\roster
\item $V$ is an irreducible $L(\lambda _n)$--module;
\item $\overline W V = 0\ $.
\endroster
\endproclaim

\proclaim
{Corollary 4}
Let $n \in \Bbb N$ and $\langle \lambda,c \rangle = n-\tfrac 3 2$. We have 
$$\overline W M(\lambda) = 
\left  \{ \aligned & M^1(\lambda ) \qquad
\text {for all}\  \lambda \in \ S_1^n \cup S_2^n, \\
& M(\lambda) \qquad \quad \text{otherwise}. \endaligned \right.$$
\endproclaim

\enddocument
\end